\documentclass[12pt]{article}
\usepackage{amssymb}
\def \C{{\mathbb C}}

\def \Z{{\mathbb Z}}
\def \p{{\mathbb P}}
\def \Q{{\mathbb Q}}

\begin{document}

\title{A computation of invariants of a rational self-map}
\author{Ekaterina Amerik}
\date{}
\maketitle
Let $V$ be a smooth cubic in $\p^5$ and let $X={\cal F}(V)$ be the variety
of lines on $V$. Thus $X$ is a smooth four-dimensional subvariety
of the grassmanian $G(1,5)$, more precisely, the zero locus of a section
of $S^3U^*$, where $U$ is the tautological rank-two bundle over $G(1,5)$. It is 
immediate from this description that the canonical class of $X$ is
trivial. Let ${\cal F}\subset V\times X$ be the universal family
of lines on $V$, and let $p:{\cal F}\to V$, $q:{\cal F}\to X$ be the 
projections. Beauville and Donagi prove in \cite{BD} that the 
Abel-Jacobi map $AJ: q_*p^*: H^4(V,\Z)\to H^2(X,\Z)$ is an isomorphism,
at least after tensoring up with $\Q$; since this is also a morphism of 
Hodge structures,
we obtain, by the Noether-Lefschetz theorem for $V$, that $Pic(X)=\Z$ for
a sufficiently general $X$.

Lines on cubics have been studied in \cite{CG}. It is shown there that
the normal bundle of
a general line on a smooth 4-dimensional cubic is
${\cal O}_{\p^1}\oplus {\cal O}_{\p^1}\oplus {\cal O}_{\p^1}(1)$ ($l$ is then 
called "a line of the first kind"),
and that some special lines  ("lines of the second kind") have normal
bundle ${\cal O}_{\p^1}(-1)\oplus {\cal O}_{\p^1}(1)\oplus {\cal O}_{\p^1}(1)$. The lines
of the second kind form a two-parameter subfamily $S\subset X$, and all
lines on $V$ are either of the first or of the second kind.

An alternative description of lines on a smooth cubic in $\p^n$ from \cite{CG}
is as follows: let $F(x_0, \dots x_n)=0$ define a smooth cubic $V$ in $\p^n$ 
and 
consider the Gauss map 
$$D_V:\p^n\to (\p^n)^*: \ x\mapsto (\frac{\partial F}{\partial
  x_0}(x): 
\dots : \frac{\partial F}{\partial x_n}(x))$$ 
(so $D_V(x)$ is the tangent hyperplane to $V$ at $x$). For a line 
$l\subset V$ there are only two possibilities: either $D_V|_l$ maps $l$ 
bijectively onto a plane conic, or $D_V|_l$ is two-to-one onto a line in 
$(\p^n)^*$. In the first case, $l$ is of the first kind and the intersection
of the hyperplanes $\cap_{x\in l} D_V(x)\subset \p^n$ is a subspace of dimension
$n-3$, tangent to $V$ along $l$; in the other case, $l$ is of the second kind
and the linear subspace $Q_l=\cap_{x\in l} D_V(x)\subset \p^n$, 
which is of course still tangent to 
$V$ along $l$,  is of greater dimension $n-2$. 
 
As remarked in \cite{V}, $X$ always admits a rational self-map.
This map is constructed as follows: from the description of the normal bundle,
one sees that for a general line $l\subset V$, there is a unique plane 
$P_l\subset \p^5$ which is tangent to $V$ along $l$. So $P_l\cap V$ is the
union of $l$ and another line $l'$, and one defines $f:X\dasharrow X$
by sending $l$ to $l'$. This is only a rational map: indeed, if $l$
is of the second kind, we have a one-parameter family, more precisely, a 
pencil of planes $\{P:\ l\subset P\subset Q_l\}$
tangent to $V$ along $l$. So the surface $S$ is the indeterminacy locus of
the map $f$. Blowing $S$ up, one obtains a resolution of singularities of $f$; 
it introduces the
exceptional divisor $E$ which is a $\p^1$-bundle over $S$, each $\p^1$
corresponding to the pencil of planes tangent  to $V$ along $l$. That
is, the "image" by $f$ of every exceptional point is a (possibly singular) 
rational curve on $X$.  

In \cite{V}, it is shown that the degree of $f$ is 16. The proof is very 
short: the space $h^{2,0}(X)$ is generated by 
a nowhere-degenerate form $\sigma$
(\cite{BD}); C.Voisin then uses a Mumford-style argument (saying that 
a family of 
rationally equivalent cycles induces the zero map on certain
differential forms) to show that
$f^*\sigma=-2\sigma$. As $(\sigma\overline{\sigma})^2$ is a volume form,
the result follows.

In holomorphic dynamics, one considers the  {\it dynamical degrees}
of 
a rational self-map (the definition below seems to appear 
first in \cite{RS}). Those are defined as follows:
let $X$ be a K\"ahler manifold of dimension $k$.
Fix a K\"ahler class $\omega$. Define 
$$\delta_l(f,\omega) = [f^*\omega^l]\cdot[\omega^{k-l}]$$ 
(the product of cohomology classes on $X$), and the $l$th
{\it dynamical degree} $$\lambda_l(f)= \lim \sup (\delta_l(f^n,\omega))^{1/n}$$
(this does not depend on $\omega$). It turns out that this is the same as
$$\rho_l(f)= \lim \sup (r_l(f^n))^{1/n},$$
where $r_l(f^n)$ is the spectral radius of the action of $f^n$ on 
$H^{l,l}(X)$. If $X$ is projective, one can restrict oneself to the
subspace of algebraic classes, replacing $\rho_l$ and $r_l$ by the 
correspondent spectral radii $r_l^{alg}$ and $\rho_l^{alg}$.
Moreover, Dinh and Sibony \cite{DS} recently proved that the $\lim
\sup$ in the definition of $\lambda_l(f)$ is actually a limit,
and also the invariance of the dynamical degrees by birational
conjugation.

The map $f$ is said to be {\it cohomologically hyperbolic} if
there is a dynamical degree which is strictly greater then all
the others. In fact the dynamical degrees satisfy a convexity
property which says that if $\lambda_l(f)$ is such a maximal
degree,
then $\lambda_i$ grows together with $i$ until $i$ reaches $l$, and
decreases thereafter. It is conjectured that the cohomological
hyperbolicity implies certain important equidistribution properties
for $f$, and this is proved in the works of Briend-Duval \cite{BrDv}
and Guedj \cite{G} when the maximal dynamical degree is the topological
degree $\lambda_{k}$.

In general it is very difficult to compute the dynamical degrees
of a rational self-map, because for rational maps it is not always
true that $(f^n)^*=(f^*)^n$ (where $f^*$ denotes the transformation induced 
by $f$ in
the cohomology.). In fact few examples of such computations are known.

If $(f^n)^*=(f^*)^n$, the map $f$ is called 
{\it algebraically stable}. The computation of its dynamical degrees
is then much easier: one only has to study the linear map $f^*$, 
and it is not necessary to look at the iterations
of $f$.

The purpose of this note is to show that the example $f$ considered in
the
beginning is algebraically stable and to compute its dynamical
degrees. The map $f$ turns out to be cohomologically hyperbolic, and the
dominating dynamical degree is in this case $\lambda_2$.

I am very much indebted to A. Kuznetsov and C. Voisin for crucially
helpful
discussions.

\

{\bf 1. Algebraic stability of the self-map}

\ 

From now on let $X={\cal F}(V)$ and $f$ be as described in the beginning.
Let $S$ be the indeterminacy locus $I(f)$, that is, the surface of
lines of the second kind.

We assume that $V$ is {\it sufficiently general}.  
The first goal is to prove that $f$ is algebraically stable.

\medskip

{\bf Lemma 1:} {\it For a general $V$, the surface $S$ is smooth.}

\medskip

{\it Proof:} Let $T$ be the parameter space of cubics in $\p^5$
Consider the incidence variety ${\cal I}\subset G(1,5)\times T$:
${\cal I}=\{(l,V): l\ is\ of\ the\ 2nd\ kind\ on\ V\}.$ It is enough to
show that ${\cal I}$ is smooth in codimension two. By homogeneity,
it is enough to verify that the fiber $I_l=pr_1^{-1}l$ is smooth in
codimension two for some $l\in G(1,5)$. Choose coordinates 
$(x_0:\dots :x_5)$ on $\p^5$ so that
$l$ is given by $x_2=x_3=x_4=x_5=0$; then $l\subset V$ means that
four coefficients (those of $x_0^3$, $x_0^2x_1$, $x_0x_1^2$ and $x_1^3$)
of the equation $F=0$ of $V$ vanish, and $l$ being of the second type means that
the 4 quadratic forms $\partial F/\partial x_i|_l$, $i=2,3,4,5$ in two 
variables
$x_0, x_1$, span a linear space of dimension 2, rather than 3. This in turn 
means that a $4\times 3$ matrix $M$, whose coefficients are (different)
coordinates on $T'\subset T$, the projective space of cubics containing $l$,
has rang 2. It follows from general facts about determinantal varieties
that $Sing(I_l)$ is the locus where $rg(M)\leq 1$, and it has codimension
four in $I_l$ (in fact, $Sing(I_l)$ is also contained in the locus of
singular cubics, \cite{CG}).

\medskip

Let us now look at the map $f: X\dasharrow X$. The first thing to 
observe is as follows:
$K_X=0$ means that $f$ cannot contract a divisor.
Indeed, let $\pi: Y\to X,\ g:Y\to X$ be the resolution of the indeterminacies
of $f$. Then the ramification divisor (i.e. the vanishing locus of the Jacobian
determinant) of $\pi$ is equal to that of $g$, so $g$
can only contract a divisor which is already
$\pi$-exceptional. By the same token, any subset contracted by $g$ must
lie in the exceptional divisor of $\pi$, in other words, anything contracted
by $f$ is in $S$.

\medskip

{\bf Lemma 2:} {\it The map $g$ does not contract a surface (neither to a
point nor to a curve), unless possibly some surface already contracted
by $\pi$. }

\medskip

{\it Proof (C. Voisin):} If it does, this surface $Z$ projects onto $S$.
Recall that $X$ is holomorphic symplectic, with symplectic form
$\sigma=AJ(\eta),$ where $\eta$ generates $H^{3,1}$ of the cubic. 
Notice that $S$ is not lagrangian: $\sigma|_S\neq 0$. This is because,
as we shall see in the next section, the cohomology class of $S$
is $5(H^2-\Delta)$, where $H$ is the hyperplane section class for the 
Pl\"ucker embedding, and $\Delta$ is the Schubert cycle of lines 
contained in a hyperplane. We have $[\sigma]\cdot\Delta=0$ in the
cohomologies, by projection formula (because $[\eta]$ is a primitive
cohomology class). But $[\sigma]\cdot H^2$ is non-zero, because 
$[\sigma\overline{\sigma}]\cdot H^2$ is the Bogomolov square of
an ample $H$. Thus so
is $\sigma|_S$.

Now recall (\cite{V}) that $f^*\sigma=-2\sigma$, that is, 
$g^*\sigma=-2\pi^*\sigma$. This is a contradiction, since
$g^*\sigma|_Z=0$, whereas $\pi^*\sigma|_Z\neq 0$.

\medskip

{\bf Remark:} The same argument with $\sigma$ shows that $g$ cannot contract
the $\pi$-exceptional divisor (even not to a surface). A somewhat more elaborate
version of it shows that $g$
does not contract any surface at all; but we won't need this to prove the
algebraic stability. See also the remark after Theorem 3.

\medskip

{\bf Theorem 3:} {\it The self-map $f$ of $X={\cal F}(V)$ described above is
algebraically stable, that is, we have $(f^n)^*=(f^*)^n$, where the upper star
denotes the action 
on the cohomologies, and $n$ is a positive integer.}

\medskip

{\it Proof:} The transformation induced by $f$ on the cohomologies is 
given by the
class of the graph $\Gamma(f)\in H^8(X\times X)$. In general,
 $(hf)^*=h^*f^*$ is implied  by $\Gamma(f)\circ \Gamma(h)$ having only one 
4-dimensional component (which is then, of course, $\Gamma(hf)$).
Recall that the correspondence $\Gamma(f)\circ \Gamma(h)$ is
$$\{(x,z)| \exists y:\ (x,y)\in \Gamma(f),\ (y, z)\in \Gamma(h)\}.$$
Let
$h=f^n$.
Suppose that $\Gamma(f)\circ \Gamma(h)$ has an extra 4-dimensional component 
$M$.
There are three possibilities for its image under the first projection:

1) It is a divisor $D$. Then $f$ is defined at a general point $d\in D$;
but fibers of $M$ over $D$ must be at least one-dimensional. That is,
$D$ is contracted by $f$ (to the indeterminacy locus of $h$). This is impossible.

2) It is a surface $Z$; fibers of $M$ over $Z$ must be two-dimensional.
But if $Z\not\subset I(f)$, there is only one $y$ corresponding to generic $x\in Z$, and
furthermore, since $Z$ is not contracted by $f$, 
$h(y)$ is at most a curve. Indeed, the locus 
$I_2(h)=\{y: dim(h(y))\geq 2\}$ is of codimension at least three in $X$. Thus a fiber
of $M$ over a generic $x\in Z$ cannot be two-dimensional. 
If $Z\subset I(f)$, one concludes in a similar way: for any $x\in Z$, the
corresponding $y$ are parametrized by a rational curve  $D_x$; for a 
general $x$, the curve  $D_x$ is not the indeterminacy locus of $h$, and,
moreover, since no surface is contracted, $D_x$ intersects $I(h)$ in a finite
number of points, all outside $I_2(h)$.
Thus again the general fiber of 
$\Gamma(f)\circ \Gamma(h)$  over $T$ 
is at most one-dimensional.

3) It is a curve $C$. By a simple dimension count as in the previous case,
we see that then $C$ must be contracted by $f$ to a point $y$ such that
$dim(h(y))=3$. Recall that $h=f^n$. We may assume that $n$ is the smallest
number with the property $dim(h(y))=3$. As for any $x\in X$, $f(x)$ is 
either a point or a rational curve, $S$ must be contained in
$f^{n-1}(y)$ as a component. This implies that $S$ is covered by rational 
curves. But, as we shall see in the next section, $K_S$ is ample; so
this case is also impossible. 

\medskip

{\bf Remark:} In fact one can show that $g$ is finite (and this obviously
implies the algebraic stability). Indeed, by the argument of Lemma 2, if 
$g(C)$ is a point, then $\sigma|_S$ vanishes on $\pi(C)$ and so 
$\pi(C)$ is a component of a canonical divisor on $S$. On the other hand,
all lines on $V$ corresponding to the points of $\pi(C)$ intersect the line 
corresponding to the point $g(C)$. This means that $\pi(C)$ is
contained in a hyperplane section of $S$. But it turns out (see the 
remark in the end of next section) that $K_S$ is even ampler than the
hyperplane section class $H_S$. To get a contradiction,
it remains to show that the vanishing locus of $\sigma|_S$ is smooth and
irreducible for generic $V$. C. Voisin indicates that this is not very
difficult to calculate explicitely using the Griffiths' description of
the primitive cohomology of hypersurfaces (here cubics) in terms of residues.

\

{\bf 2. Some cohomology classes on $X$}

\

Recall that $X\subset G(1,5)$ is the zero-set of a section of the
bundle
$S^3U^*$ (of rang 4). 
We shall describe the action of $f^*$ on those cohomology classes
which are restrictions of classes on $G(1,5)$.

So let $H=c_1(U^*)|_X=[{\cal O}_{G(1,5)}(1)]|_X$; 
$\Delta=c_2(U^*)|_X$ (so $\Delta$ is the class of the subset
of points $x\in X$ such that the corresponding lines $l_x\subset \p^5$
lie in a given hyperplane). 

\medskip

{\bf Lemma 4:} $H^4=108;\ H^2\Delta=45;\  \Delta^2=27.$

\medskip

{\it Proof:} $[X]=c_4(S^3U^*)$ and   
$$c_4(S^3U^*)=9c_2(U^*)(2c_1^2(U^*)+c_2(U^*)).$$
Let $\sigma_i=c_i(U^*)$; those
are Schubert classes on $G(1,5)$: 
$\sigma_1=\{$lines intersecting a given $\p^3 \}$, 
$\sigma_2=\{$lines lying in a given hyperplane$\}$.

Notice that $c_1(U^*)$ is the hyperplane section class in the Pl\"ucker
imbedding and restricts as such onto a "sub-grassmanian" of lines lying
in a linear subspace. So
$$H^4=9\sigma_2\sigma_1^4(2\sigma_1^2+\sigma_2)=
18 deg(G(1,4))+9 deg(G(1,3))=108.$$ 
The other equalities are proved similarly..

\medskip

For generic $X$, the class $H$ generates $H^{1,1}_{alg}$  over $\Z$,
and $H^3$ generates $H^{3,3}_{alg}$ over $\Q$.

One can also show that $H^2$ and $\Delta$ generate $H^{1,1}_{alg}$ over $\Q$;
we do not need this for our computation, but see the remark in
the end of the last section.

\

The next step is to compute the indeterminacy surface $S$.

Let $A$ denote the underlying vector space for $\p^5\supset V$, and let 
$L\subset A$ be the underlying vector space for a line $l\subset V$.
Write the Gauss map associated to $F$ defining $V$ as
$D_F:A\to S^2A^*$. The condition $l\subset V$ implies
that this map descends to $\phi: A/L\to S^2L^*$:
$$
\begin{array}{ccccc} 
A & \rightarrow & S^2A^*\\
\downarrow && \downarrow \\
A/L\ & \rightarrow & S^2L^*
\end{array}
$$

For a general $L$, the map $\phi$ is surjective with one-dimensional 
kernel; and it has two-dimensional kernel exactly when $L$ is the
underlying vector space of a line of the second kind. Now $L$ is the fiber
over a point $x_l\in X$
of the restriction $U_X$ of $U$ to $X$; so, globalizing and 
dualizing, we get the following resolution for $S$:
$$0\to S^2U_X\to Q^*_X\to M\to M\otimes {\cal O}_S\to 0,$$
where $Q^*_X$ is the restriction to $X$ of the universal quotient bundle
over $G(1,5)$, and $M$ is a line bundle. In fact
$$[M]=det(Q^*_X)-det(S^2U_X)=-H +3H=2H.$$ 
So one has a resolution for the ideal sheaf ${\cal I}_S$:
$$0\to S^2U_X(-2H)\to Q^*_X(-2H)\to {\cal I}_S\to 0.$$

The cohomology class $[S]$ can be computed by the Thom-Porteous
formula: a partial case of this formula identifies the cohomology class of 
the degeneracy locus $D_{e-1}(\phi)$,
where a vector bundle morphism $\phi: E\to F$ ($rg(E)=e\leq rg(F)=f$)
is not of maximal rang, as $c_{f-e+1}(F-E)$, where we put formally
$$c(F-E)=1+c_1(F-E)+c_2(F-E)+\dots = c(F)/c(E).$$
So $$[S]=c_2(Q^*_X-S^2U_X)=c_2(Q^*_X)-c_2(S^2U_X)+2H c_1(S^2U_X)=$$
$$=H^2-\Delta-4\Delta-2H^2+6H^2=5(H^2-\Delta).$$

One can also make a direct calculation using the equality 
$$[S]=-c_2({\cal O}_S)=c_2({\cal I}_S)$$
(see e.g. \cite{F},
chapter 15.3).
Later, we shall need its extension
$$i_*c_1(N_{S,X})=i_*K_S=c_3({\cal I}_S)$$
where $i:S\to X$ is the usual embedding. From the resolution of the ideal
sheaf, one computes
$$c_3({\cal I}_S)=c_3(Q^*_X(-2H))-c_3(S^2U_X(-2H))-
c_2({\cal I}_S)c_1(S^2U_X(-2H))=20H^3-27H\Delta.$$

Thus we have obtained the following

\medskip

{\bf Proposition 5:} {\it The Chern classes of ${\cal I}_S$ are:
$c_1({\cal I}_S)=0$; $c_2({\cal I}_S)=5(H^2-\Delta);$ 
$c_3({\cal I}_S)=
20H^3-27H\Delta = \frac{35}{4}H^3$.}

\medskip

The last equality is obtained by taking the intersection with $H$
and using lemma 4.

\medskip

{\bf Remark:} One can also remark that $S$ being the degeneration locus
of a map $\phi:S^2U\to Q^*$, its normal bundle is isomorphic to 
$(Ker(\phi|_S))^*\otimes (Coker(\phi|_S))$, so its determinant is
$c_1(Coker(\phi|_S))-2c_1(Ker(\phi|_S))=c_1(Q^*|_S)-c_1(S^2U|_S)
-c_1(Ker(\phi|_S))=2H_S-c_1(Ker(\phi|_S)),$
in particular it is at least as positive than $2H_S$. It is easy to see
that the cohomology class of the intersection $HS$ is
$\frac{35}{12}H^3$; this suggests that the canonical class of $S$ is
 equal to $3H_S$, but I have not checked this (one probably can write an
analogue of the Thom-Porteous formula, but I could not find a 
reference, and in 
our particular case the direct calculation is easy).

\

{\bf 3. Computation of the inverse images}

\

Let $\pi: \tilde{X}\to {X}$ be the blow-up of $S$; this resolves
the singularities of $f$, that is, $g=f\circ\pi: \tilde{X}\to {X}$
is holomorphic. Let $E$ denote the exceptional
divisor of $\pi$. We have $g^*H=a\pi^*H-bE$; the number $a$ is, by
algebraic stability, 
the first
dynamical degree we are looking for. In fact the numbers $a$ and $b$ can be
computed by a simple geometric argument; we begin by sketching this 
computation.

\medskip

{\bf Proposition 6} $a=7$ and $b=3$.

\medskip

{\it Proof:} Consider the two-dimensional linear sections $V\cup \p^3$ of $V$; 
they form
a family of dimension $dim(G(3,5))=8$. A standard dimension count shows
that in this family, there is a one-parameter subfamily of cones over
a plane cubic (indeed, one computes that out of $\infty^{19}$ cubics in
$\p^3$, $\infty^{12}$ are such cones, so being a cone imposes only 7
conditions on a cubic surface in $\p^3$). Furthermore, some (a finite number)
of those plane cubics are degenerate, so on $V$ we have a finite number of 
cones over a plane cubic with a double point (as it was already mentioned in
the first section). Thus on $X$, we have a one-parameter family of plane cubic
curves $C_t$.

Another way to see this is to remark that the lines passing through
a point of $V$ sweep out a surface which is a cone over the intersection 
of a cubic and a quadric in $\p^3$ (\cite{CG}). The quadric sometimes
degenerates in the union of two planes, so the intersection becomes
a union of two plane cubics. Moreover, the two plane cubics have
three points in common (e.g. from the arithmetic genus count).

The resulting curves $C_t$ on $X$ are clearly $f$-invariant, and when smooth, 
map
4:1 to themselves (indeed, through a given point $p$ of an elliptic
curve
$C\in \p^2$, there are 4 tangents
to $C$, because the projection from $p$ is a 2:1 map from $C$ to $\p^1$ and
so has 4 ramification points). A smooth $C_t$ has 3 points in common with $S$.
Let $C_0$ be singular; then it has at least one point $x$ in common with $S$,
and is the image by $g$ of the exceptional $E_x\cong \p^1\subset \tilde{X}$ 
over 
this point; $g$ is 1:1 from $E_x$ to $C_0$. 

It follows that $E_xg^*H=3$, so $b=3$ (recall that $E$ restricts to itself as
${\cal O}_E(-1)$, hence $E_xE=-1$). Also, $C_tH=C_tE=3$ and $C_tg^*H=12$,
from where $a=7$.  

\medskip

However, I could not find such a simple way to compute the 
higher dynamical degrees. The following approach is suggested by 
A. Kuznetsov.

On $\tilde{X}$, one has a vector bundle $P$ of rank three:
its fiber over $x_l$ is the underlying vector space of the plane tangent 
to $V$ along $l$ when $l$ is of the first kind, and we take the obvious
extension
to the exceptional divisor. One has $\pi^*U_X\subset P$.

First of all, let us describe the quotient: let $E$ denote the exceptional
divisor of $\pi$. Then one can see that we obtain  $P$ as the
following extension:
$$
\begin{array}{ccccccccc}
0 & \to& \pi^*U_X & \to & P & \to & M^*(E)& \to & 0 \\
&& || && \downarrow && \downarrow && \\
0 & \to& \pi^*U_X & \to & A\otimes {\cal O}& \to & \pi^*Q_X & \to & 0,
\end{array}
$$
where the map $M^*(E)\to \pi^*Q_X$ is the dual of the pull-back of 
the map $Q^*_X\to M\otimes {\cal I}_S$ from the resolution of $S$ in the
last section.

(Alternatively, one can compute the cokernel of the natural inclusion
$\pi^*U_X  \to  P$ using the exact sequences from the proof of the 
next proposition
and the knowledge of $c_1(g^*U^*_X)=7\pi^*H-3E$).

\medskip

{\bf Proposition 7:} {\it The bundle $g^*U^*_X$ fits into an exact sequence}
$$0\to (M^*)^{\otimes 2}(2E)\to P^*\to g^*U^*_X\to 0.$$

\medskip
 
{\it Proof:} Consider the projective bundle $p:\p_{\tilde{X}}(P)\to \tilde{X}$.
It is equipped with a natural map $\psi: \p_{\tilde{X}}(P)\to \p(A)$, and
the inverse image of ${\cal O}_{\p(A)}(1)$ is the tautological line bundle 
${\cal O}_{\p_{\tilde{X}}(P)}(1)$, meaning that its direct image under the 
projection
to $\tilde{X}$ is $P^*$. Let $h=c_1({\cal O}_{\p_{\tilde{X}}(P)}(1)).$
We have a natural inclusion $\p_{\tilde{X}}(\pi^*U_X)\subset\p_{\tilde{X}}(P)$
.
Denote by $d$ the class of $\p_{\tilde{X}}(\pi^*U_X)$ considered as a divisor
on $\p_{\tilde{X}}(P)$. The subvariety of $\p_{\tilde{X}}(P)$ swept out
by the lines $f(l)$ (where $l$ is a fiber of $\p_{\tilde{X}}(\pi^*U_X)$),
is a divisor; let us denote it by $\Sigma$. We clearly have $2d+[\Sigma]=3h$.
Furthermore, $\Sigma$ itself is a projective bundle over
$\tilde{X}$, namely, it is the projectivization of $g^*U_X$,
in the sense that the vector bundle $g^*U^*_X$ is the direct image of 
${\cal O}_{\Sigma}(1)=
{\cal O}_{\p_{\tilde{X}}(P)}(1)|_{\Sigma}$ under the projection to $\tilde{X}$.  
On $\p_{\tilde{X}}(P)$, we have the exact sequence
$$0\to {\cal O}_{\p_{\tilde{X}}(P)}(h-[\Sigma])\to {\cal O}_{\p_{\tilde{X}}(P)}(1)
\to {\cal O}_{\Sigma}(1)\to 0.$$
Also, $h-[\Sigma]=2(d-h)$. So our proposition follows by taking the direct
image of the above exact sequence, once we show that 
$${\cal O}_{\p_{\tilde{X}}(P)}(h-d)=p^*M(-E).$$ 
As for this last statement, it is clear that 
${\cal O}_{\p_{\tilde{X}}(P)}(h-d)=p^*F$ for some line bundle $F$
on $\tilde{X}$, because it is trivial on the fibers. We deduce that
$F=M(-E)$ from the exact sequence
$$0\to {\cal O}_{\p_{\tilde{X}}(P)}(h-d)\to {\cal O}_{\p_{\tilde{X}}(P)}(1)
\to {\cal O}_{\p_{\tilde{X}}(P)}(1)|_{\p_{\tilde{X}}(\pi^*U_X)}\to 0:$$
indeed, it pushes down to 
$$0\to F\to P^*\to \pi^*U^*_X\to 0.$$

We are ready to prove the main result:

\medskip

{\bf Theorem 8} {\it We have $f^*H=7H$, $f^*H^2=4H^2+45\Delta$, 
$f^*\Delta=31\Delta$ and $f^*H^3=28H^3$. The dynamical degrees of $f$ are
therefore 7, 31, 28 and 16. In particular, $f$ is cohomologically 
hyperbolic.}

\medskip

{\it Proof:} We have already seen (proposition 6) that $f^*H=\pi_*g^*H=7H$.
Notice that $\pi_*(\pi^*HE)=0$ and $\pi_*E^2=-S$ (this is because 
${\cal O}(E)$ restricts
to $E$ as ${\cal O}_E(-1)$, where we view $E$ as a projective bundle
over $S$). This gives 
$$f^*H^2=49H^2-9S=49H^2-45(H^2-\Delta)=4H^2+45\Delta.$$
Further, $g^*\Delta=c_2(g^*U_X)$, and we can find the latter from
the exact sequences describing $g^*U_X$ and $P$.
From $$0\to \pi^*U_X \to P\to M^*(E)\to 0,$$ we get
$c_2(P)=\pi^*\Delta-\pi^*H(-2\pi^*H+E)=\pi^*\Delta-\pi^*HE+2\pi^*H^2,$ and from 
$$0\to g^*U_X \to P\to M^{\otimes 2}(-2E)\to 0,$$ 
$$c_2(g^*U_X)=c_2(P)+(4\pi^*H-2E)(7\pi^*H-3E)=\pi^*\Delta+30\pi^*H^2+6E^2-27\pi^*HE,
$$
that is, $$f^*\Delta=\Delta+30H^2-30(H^2-\Delta)=31\Delta.$$
Finally, 
$$g^*H^3=(7\pi^*H-3E)^3=343\pi^*H^3-441\pi^*H^2E+189\pi^*HE^2-27E^3.$$
We have $\pi_*(\pi^*H^2E)=0$ and 
$$\pi_*(\pi^*HE^2)=-5H(H^2-\Delta)=-\frac{35}{12}H^3.$$
Let us compute $\pi_*E^3$: this is $\pi_*\xi^2$, where $\xi$ is the tautological
class on $E$ viewed as a projective bundle $\p_S(N)$ over $S$. Let $r$ be the
projection of $E$ to $S$, i.e. the restriction of $\pi$ to $E$. We have
$\xi^2+c_1(r^*N)\xi+c_2(r^*N)=0,$ which yields $\pi_*E^3=-i_*c_1(N)$, $i$ being
the embedding of $S$ into $X$. But the latter is just $c_3(i_*{\cal O}_S)$
(\cite{F}). We have computed this class in the last section: 
$$\pi_*E^3=27H\Delta-20H^3=-\frac{35}{4}H^3.$$
Putting all this together, we get
$$f^*H^3=343H^3-\frac{35}{4}(63-27)H^3=28H^3.$$

This finishes the proof of Theorem 8.

\medskip

{\bf Remark 9:} Notice that the eigenvectors of $f^*$ on the invariant subspace
generated by $H^2$ and $\Delta$ are orthogonal with respect to the 
intersection form. There are reasons for this; moreover, at least 
for $X$ generic the map $f^*$ on the 
whole cohomology
group $H^4(X)$ is self-adjoint with
respect to the intersection form. One can see  it as follows: $X$ is a 
deformation 
of
the punctual Hilbert scheme $Hilb^2(S)$ of a $K^3$ surface \cite{BD}; 
this implies
a decomposition of $H^4(X)$ into an orthogonal direct sum of Hodge substructures
$$H^4(X)=<H^2>\oplus (H\cdot H^2(X)^0) \oplus S^2(H^2(X)^0),$$
where $H^2(X)^0$ denotes the orthogonal to $H$ in $H^2(X)$. For $X$ generic, 
the second
summand is an irreducible Hodge structure and the third one is a sum
of an irreducible one and a one-dimensional subspace: this follows from
the fact that $H^2(X)^0=H^4(V)^{prim}$ (where $V$ is the corresponding cubic 
fourfold and "prim" is primitive cohomologies) and a theorem of Deligne
which says that the closure of the monodromy group (of a general 
Lefschetz pencil) acting
on $H^4(V)^{prim}$ is the full orthogonal group. Therefore the decomposition
rewrites as 
$$H^4(X)=<H^2, \Delta>\oplus H\cdot H^2(X)^0 \oplus V,$$
where the last two summands are simple, and so $f^*$ must be a homothety
on each of them.

It also follows from this decomposition that on $X$ generic,
the space of algebraic cycles of codimension two
is only two-dimensional, and that all lagrangian surfaces must have the 
cohomology class proportional to $\Delta$. In particular, since 
$f^*$ and $f_*$ must preserve the property of being lagrangian,
this explains why $\Delta$ is an eigenvector of $f^*$ and $f_*$.

{}

\end{document}